\documentclass{amsart}
\usepackage[utf8]{inputenc}
\usepackage{amsmath}
\usepackage{latexsym,amsthm,amsmath,amssymb,amsfonts,epsfig,psfrag}
\usepackage{mathrsfs}
\usepackage{caption}

\usepackage[dvipsnames]{xcolor}

\usepackage[backref=false, colorlinks, linkcolor=blue, citecolor=Blue]{hyperref}
\usepackage{cleveref}
\usepackage{nicefrac}
\usepackage[inline]{enumitem}
\usepackage{enumitem}

\usepackage{textcomp}

\usepackage{tikz}
\definecolor{SolutionColor}{gray}{0.85}

\usepackage{euscript}

\usepackage{multirow}
\usepackage{graphicx}

\usepackage{array,float}

\usepackage{tikz}
\usetikzlibrary{shapes,arrows}
\usetikzlibrary{fit,positioning}
\usetikzlibrary{patterns,decorations.pathreplacing}
\usepackage{xkeyval}
\usepackage{moreverb}
\usepackage{pgf}
\usepackage{epic}
\usepgfmodule{shapes,plot,decorations}

\usepackage{booktabs} 
\usepackage[normalem]{ulem}

\newcommand{\Q}{\mathbb{Q}}
\newcommand{\R}{\mathbb{R}}

\newcommand{\HH}{\mathbb{H}}

\newcommand{\G}{\mathbf{G}}

\usepackage[warn]{mathtext}

\newcommand{\OOO}{\mathcal{O}}

\newcommand{\id}{\mathrm{Id}}

\newcommand{\Isom}{\mathrm{Isom}}
\newcommand{\Cyc}{\mathrm{Cyc}}

\newcommand{\Ad}{\mathrm{Ad}}

\newcommand{\tr}{\mathrm{tr}}

\newtheorem{theorem}{Theorem}[section]

\theoremstyle{remark}
\newtheorem{remark}[theorem]{Remark}

\theoremstyle{remark}

\renewcommand{\qed}{\hfill$\scriptstyle\blacksquare$}

\title{Quasi-arithmetic hyperbolic Coxeter prisms}	
 	
\author{Nikolay Bogachev}
\address{Department of Computer and Mathematical Sciences, University of Toronto Scarborough, 1095 Military Trail, Toronto, ON M1C 1A3, Canada}
\address{Institute for Information Transmission Problems, Moscow, Russia}
\email{n.bogachev@utoronto.ca}
 
\author{Khusrav Yorov}
\address{Visual Computing Center, KAUST, Thuwal 23955-6900, Saudi Arabia}
\email{khusrav.yorov@kaust.edu.sa}

\begin{document}

\begin{abstract}
    In 1974, Kaplinskaja classified all simplicial straight hyperbolic Coxeter prisms. In this paper, we determine precisely which of these prisms are properly quasi-arithmetic or arithmetic. We also present some observations regarding commensurability classes and systoles of the associated orbifolds.
\end{abstract}
\maketitle

\section{Introduction}

This paper concerns arithmetic properties of {\em simplicial hyperbolic Coxeter prisms}, that is, finite-volume hyperbolic Coxeter polyhedra combinatorially equivalent to the product of a simplex and a segment. Recall that hyperbolic Coxeter polyhedra, i.e. convex polyhedra with dihedral angles of the form $\pi/m$, are canonical fundamental domains for discrete groups generated by reflections in hyperplanes of hyperbolic (or Lobachevsky) spaces~$\HH^n$. 

Any simplicial prism has precisely two simplicial facets, its {\em bases}. In contrast to the Euclidean setting, a hyperbolic Coxeter prism can have at most one base orthogonal to all  adjacent facets. In the presence of such a base, which is necessarily compact, the prism is said to be {\em straight}. Any hyperbolic Coxeter prism can be cut into two straight prisms sharing a common such base; in particular, since compact Coxeter simplices do not exist in $\HH^{\ge 5}$,  finite-volume hyperbolic Coxeter prisms cease to exist in dimensions $\ge 6$.

Compact and finite-volume straight hyperbolic Coxeter prisms were classified by Kaplinskaja \cite{Kap74} in 1974; see Table~\ref{tab:all-prisms}. 
In $\HH^3$, these prisms provide not just an infinite family of Coxeter polyhedra, but even infinitely many {\em commensurability classes} of finite-covolume reflection groups. On the other hand, in $\HH^4$ and $\HH^5$ there are only finitely many such prisms. Our goal is to precisely determine which hyperbolic reflection groups associated to straight Coxeter prisms are {\em (quasi-)arithmetic}. 

A hyperbolic lattice $\Gamma < \Isom(\HH^n) = G$ is said to be {\em quasi-arithmetic} if there is a totally real number field $\mathbf{k} \subset \R$ and a $\mathbf{k}$-group $\G$ such that $\G(\R)$ is isogenous to $G$ with $\Gamma$ virtually contained in $\G(\mathbf{k})$ via this isogeny, but $\G^{\sigma}(\R)$ is compact for every non-identity embedding $\sigma \colon \mathbf{k} \to \R$. In this case, the adjoint trace field of $\Gamma$ coincides with $\mathbf{k}$, and $\Gamma$ is
{\em arithmetic} if and only if it is commensurable to $\G(\OOO_\mathbf{k})$ where $\OOO_\mathbf{k}$ is the ring of integers of the field $\mathbf{k}$. The latter condition is equivalent to $\tr(\Ad \gamma) \in \OOO_{\mathbf{k}}$ for each $\gamma \in \Gamma$ where $\Ad$ denotes the adjoint representation of the Lie group $G$. 
We say that a lattice $\Gamma$ is {\em properly quasi-arithmetic} if it is quasi-arithmetic but not arithmetic. We will conflate reflection groups, their fundamental Coxeter polyhedra and associated {\em reflective orbifolds} (i.e. quotients of $\HH^n$ by reflection groups) when speaking about their arithmetic or group-theoretic properties.

The notion of quasi-arithmeticity was first suggested by Vinberg in 1967 in his fundamental work \cite{Vin67} on hyperbolic reflection groups where, applying his original techniques, he provided several properly quasi-arithmetic reflection groups in $\HH^n$ for $n=3,4,5$. We would like to stress that some hyperbolic Coxeter prisms were among those examples. Moreover, even earlier, in 1966, Makarov \cite{Mak66} provided the first examples of nonarithmetic lattices in $\Isom (\HH^3)$ via finite-volume noncompact Coxeter prisms. This remarkable work of Makarov motivated further study (initiated by Vinberg) of hyperbolic reflection groups and their arithmetic properties.

It is well known that there are many nonarithmetic hyperbolic lattices, including some constructions that apply in all dimensions $n \geqslant 2$. The first is due to Gromov and Piatetski--Shapiro \cite{GPS} (see also Raimbault \cite{R13} and Gelander--Levit~\cite{GL14}), who constructed nonarithmetic hyperbolic {\em hybrids} by gluing pairwise incommensurable arithmetic blocks. Later, Thomson~\cite{Thomson} observed that such hybrids are not quasi-arithmetic. Other nonarithmetic lattices were constructed by Agol~\cite{agol2006systoles}, Belolipetsky--Thomson \cite{belolipetskythomson}, and Bergeron--Haglund--Wise \cite{bergeronhaglundwise} (see also Mila~\cite{Mila18}), giving rise to finite-volume hyperbolic manifolds with arbitrarily small {\em systole}. This construction produces properly quasi-arithmetic hyperbolic lattices as was also shown by Thomson \cite{Thomson}, allowing him to distinguish these two classes of nonarithmetic lattices even up to commensurability. It is also worth mentioning that Douba \cite{Douba} recently provided another family of hyperbolic manifolds with small systoles which are not quasi-arithmetic and thus are not commensurable to those constructed by Agol, Belolipetsky--Thomson, and Bergeron--Haglund--Wise, but are created by taking hybrids of the manifolds constructed by the latter authors. 

To conclude the discussion of quasi-arithmetic hyperbolic lattices in a general setting, we mention the following two results. Emery \cite{Emery17} recently demonstrated that the covolume of any quasi-arithmetic lattice is a rational multiple of the covolume of an arithmetic lattice associated to the same ambient group. Besides this, it was shown by Belolipetsky et al. \cite[Theorem~1.7]{BBKS21} that (quasi-)arithmeticity is inherited by totally geodesic suborbifolds; more precisely, if $M$ is a quasi-arithmetic hyperbolic orbifold with adjoint trace field $\mathbf{k}(M)$, and $N \subset M$ is a finite volume totally geodesic suborbifold of dimension $m \geqslant 2$ with adjoint trace field $\mathbf{k}(N)$, then $N$ is hyperbolic and quasi-arithmetic, with $\mathbf{k}(M) \subseteq \mathbf{k}(N)$.

Returning to the more specific setting of quasi-arithmetic reflection groups, it was shown in \cite{BK21} that Coxeter faces of quasi-arithmetic Coxeter polyhedra are also quasi-arithmetic with the same ground field. Later in \cite{BD23}, the compact right-angled Löbell polyhedra $L_n \subset \HH^3$ were studied, and it was observed that $L_n$ is quasi-arithmetic if and only if $n =  5, 6, 8$, or $12$. In \cite{BGV23}, quasi-arithmetic ideal right-angled twisted antiprisms and associated hyperbolic link complements were classified.

Dotti and Kolpakov \cite{dottikolpakov2022} recently showed that there are infinitely many maximal properly quasi-arithmetic reflection groups in $\HH^2$, while it is known that there are only finitely many maximal arithmetic hyperbolic reflection groups for all dimensions $n \geqslant 2$ (see Vinberg \cite{vinberg84absence}, Esselmann \cite{Esselmann}, Agol--Belolipetsky--Storm--Whyte \cite{ABSW}, Nikulin \cite{Nikulin}). Quasi-arithmeticity was recently used in \cite{BDR23} to prove that Makarov's \cite{Mak68} polyhedra give rise to infinitely many commensurability classes of cocompact reflection groups in $\HH^4$ and $\HH^5$.

The main result of this paper is the following classification of arithmetic and properly quasi-arithmetic hyperbolic straight Coxeter prisms.

\begin{theorem}\label{th:main}
There are only $25$ properly quasi-arithmetic and $37$ arithmetic straight Coxeter prisms in hyperbolic space $\HH^{n}$, where $n \in \{3,4,5\}.$ 
\end{theorem}

\begin{remark}
In $\HH^{3}$, we have exactly $31$ arithmetic and $21$ properly quasi-arithmetic straight Coxeter prisms; see Tables \ref{tab:type1-arith}, \ref{tab:type1-quasi}, \ref{tab:type2-4-quasi}, \ref{tab:types-5-11}. In $\HH^{4}$, there are only $6$ arithmetic and $4$ properly quasi-arithmetic straight Coxeter prisms, and in $\HH^{5}$, we have $3$ arithmetic and no properly quasi-arithmetic straight Coxeter prisms; see Table~\ref{tab:other-types}. 
\end{remark}

Let us now consider the more general case when a Coxeter prism is glued from two straight prisms along a common base. Recent results of Belolipetsky et al.~\cite{BBKS21} allow us to rule out quasi-arithmeticity of the significantly many such Coxeter prisms by considering the adjoint trace field of totally geodesic subspaces. One can observe that in $\HH^3$ there are only three infinite families of compact Coxeter prisms (see Table \ref{tab:all-prisms}), which we will denote by $P^{j}_{k,l,m}$ where $j = 1, 2, 3$ denotes the type (so $j$ corresponds to the $j^{\textrm{th}}$ diagram in the list) of the prism, and $k,l,m$ are the diagram parameters. Let $\Gamma^{j}_{k,l,m}$ be the associated reflection groups, and consider the $P^{j}_{k,l,m}$ as hyperbolic reflective orbifolds $\HH^n/\Gamma^{j}_{k,l,m}$.   

\begin{theorem}\label{th:non-quasi}
    For fixed $k,l \leqslant 3$ and $m$, let $P^{j,3}_{k,l,m}$, with $j = 1, 2$, be a Coxeter prism obtained by glueing two straight Coxeter prisms $P^{j}_{k,l,m}$ and $P^{3}_{k,l,m}$ along its common $(k,l,m)$ triangular face. If $m$ is coprime with $5$, then the reflection group of $P^{j,3}_{k,l,m}$ is not quasi-arithmetic.
\end{theorem}

The next theorem is a well-known fact in some special cases, see, for instance, \cite[Section~4.7.3]{McR} where commensurability classes of type-$1$ (compact) prisms from Table~\ref{tab:all-prisms} are studied. The systolic part of this theorem is related to the paper \cite{BD23} where it was shown that the systoles of compact Löbell (hyperbolic) $3$-manifolds approach~$0$.
 
\begin{theorem}\label{th:systoles}
    Given fixed $j, k, l$, the sequence $\{\Gamma^{j}_{k,l,m}\}_{m=1}^{\infty}$ gives rise to infinitely many commensurability classes of cocompact reflection groups in $\HH^3$. Moreover, the systoles of the Coxeter orbifolds $P^{j}_{k,l,m}$ approach $0$ as $m \to +\infty$.
\end{theorem}

\noindent We would like to stress that commensurability classes of the three infinite families of noncompact Coxeter prisms were well studied by Dotti--Drewitz--Kellerhals \cite{DDK23}.

\subsection*{Organization of the paper}
Section~\ref{sec:prelim} contains preliminaries. The proof of Theorem \ref{th:main} follows in Section \ref{sec:main-proof}. The proofs of Theorems \ref{th:non-quasi} and \ref{th:systoles} are presented in Sections \ref{sec:non-quasi} and \ref{sec:systoles}, respectively.

\subsection*{Acknowledgements}
We are grateful to Sami Douba for useful remarks.
 
\begin{table} 
\def\arraystretch{1.25} 
\centering
\begin{tabular}{|c|c|}
\hline
Compact prisms in $\HH^3$ (types $1$--$4$)\label{1-4} & Noncompact prisms in $\HH^3$ (types $5$--$11$)\label{5-11}\\ 
\hline
&\\

\psfrag{k}[position = 1.0]{\scriptsize $m$}
\psfrag{k1}{\scriptsize $k$}
\psfrag{k2}[position = 1.0]{\scriptsize $l$}
\psfrag{1}{\small $2\leqslant k,l\leqslant 5$}
\psfrag{2}{\small $\frac{1}{k}+\frac{1}{l}+\frac{1}{m}<1$}
\psfrag{3}{\small $2\leqslant k,l \leqslant 3$}
\psfrag{4}{\small  $\frac{1}{k}+\frac{1}{l}+\frac{1}{m}<1$}
\psfrag{5}{\small  $2\leqslant k,l \leqslant 3$}
\psfrag{6}{\small   $\frac{1}{k}+\frac{1}{l}+\frac{1}{m}<1$}
\psfrag{7}{\small $m=4,5$}
\raisebox{1cm}{\epsfig{file=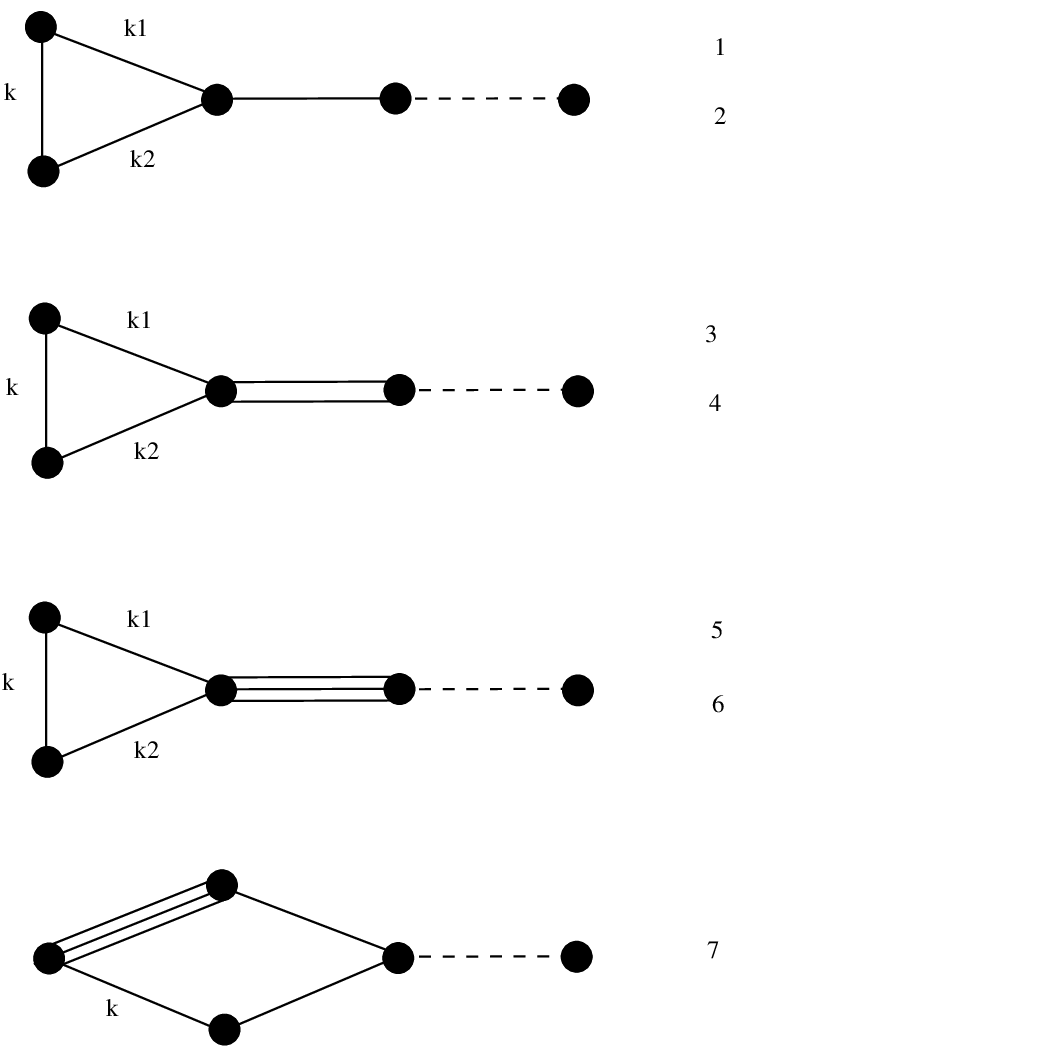,width=0.4\linewidth}}&

\psfrag{k}[position = 1.0]{\scriptsize $m$}
\psfrag{k1}{\scriptsize $k$}
\psfrag{k2}[position = 1.0]{\scriptsize $l$}
\psfrag{6}{\scriptsize $6$}
\psfrag{1}{\small $m>3$}
\psfrag{2}{\small $m>4$}
\psfrag{3}{\small $m>6$}
\psfrag{4}{\small $m>2$}
\psfrag{5}[scale = 0.5]{\small $k=3,4$}
\psfrag{7}{\small $k=3,4,5$}
\psfrag{8}[scale = 0.09]{\small $l=4,5,6$}
\epsfig{file=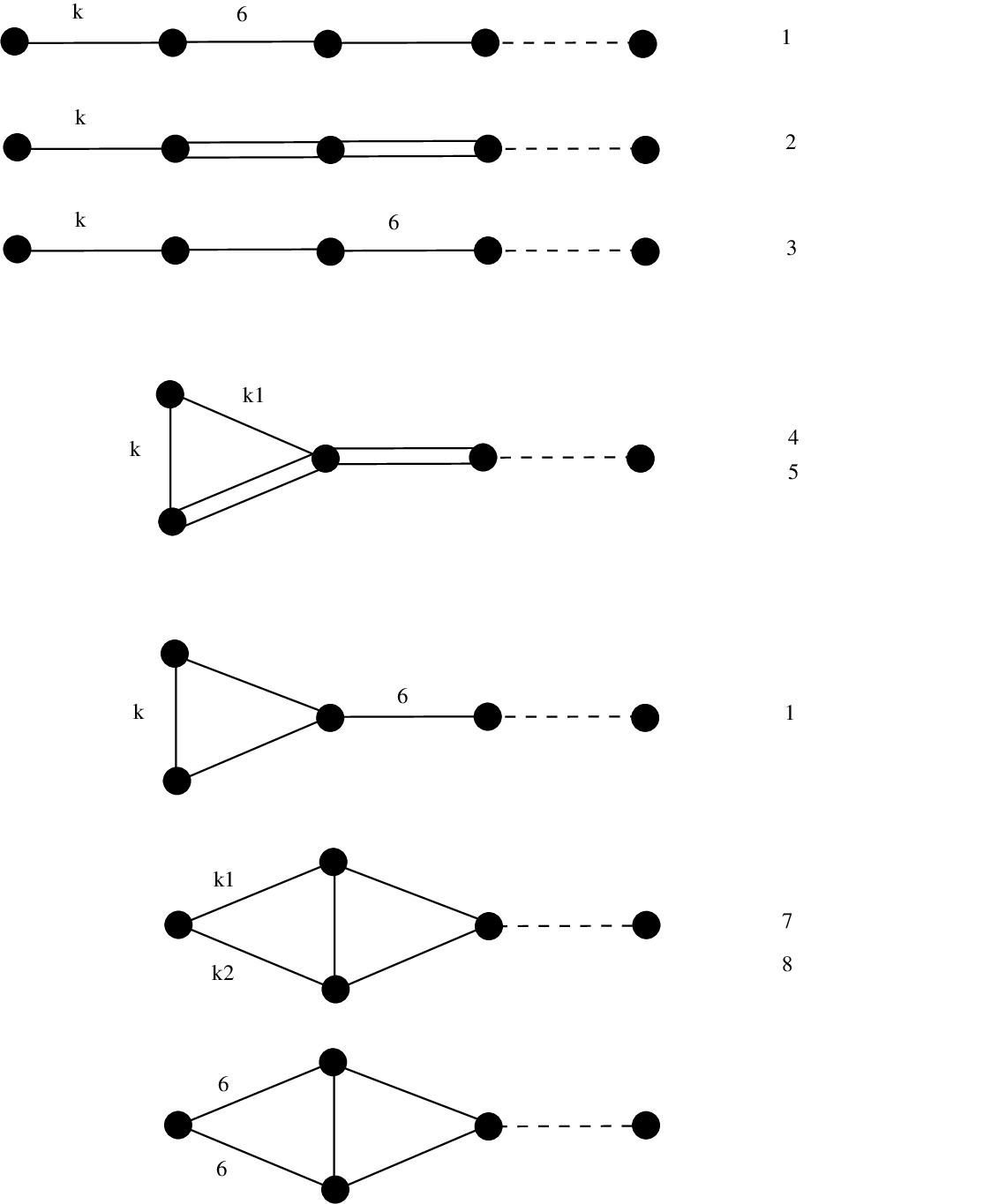,width=0.4\linewidth}\\
\hline
\hline
Compact prisms in $\HH^4$ (types $12$--$16$)\label{12-16} & Noncompact prisms in $\HH^4$ (types $17$--$20$)\label{17-20}\\ 
\hline
&\\
\psfrag{k1}{\scriptsize $k$}
\psfrag{1}{\small $k=4,5$}
\raisebox{0.5cm}{\epsfig{file=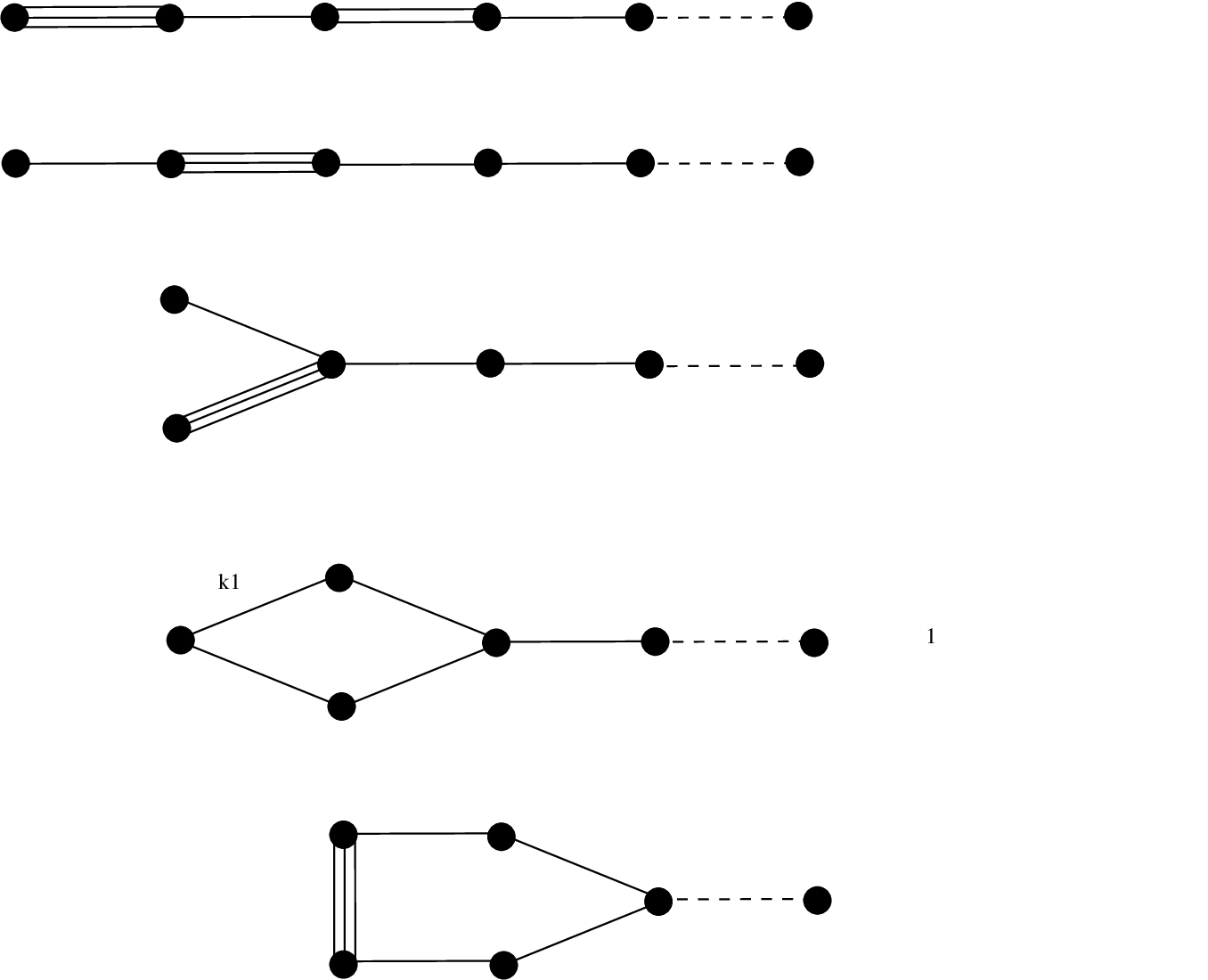,width=0.4\linewidth}}&
\psfrag{k1}{\scriptsize $k$}
\psfrag{1}{\small $k=3,4,5$}
\epsfig{file=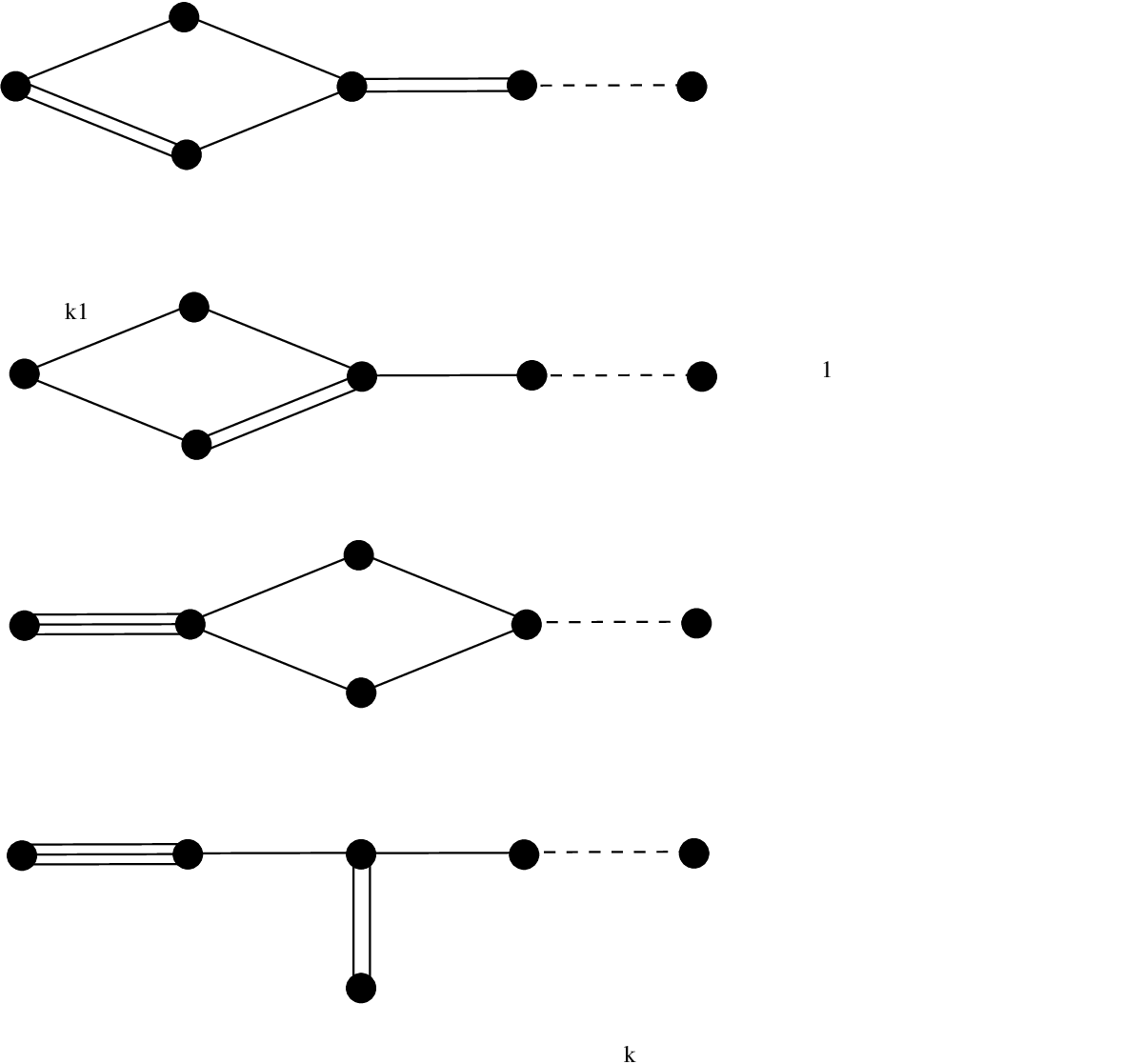,width=0.4\linewidth}\\
\hline
\hline
Compact prisms in $\HH^5$ (types $21$--$22$)\label{21-22} & Noncompact prisms in $\HH^5$ (types $23$--$24$)\label{23-24}\\ 
\hline
&\\
\psfrag{k1}{\scriptsize $k$}
\psfrag{1}{\small $k=3,4$}
\epsfig{file=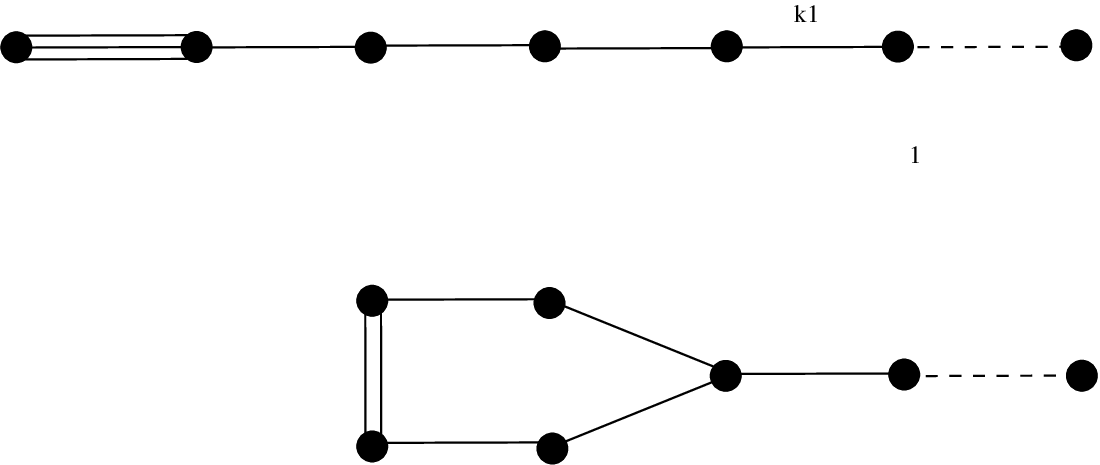,width=0.35\linewidth}&
\epsfig{file=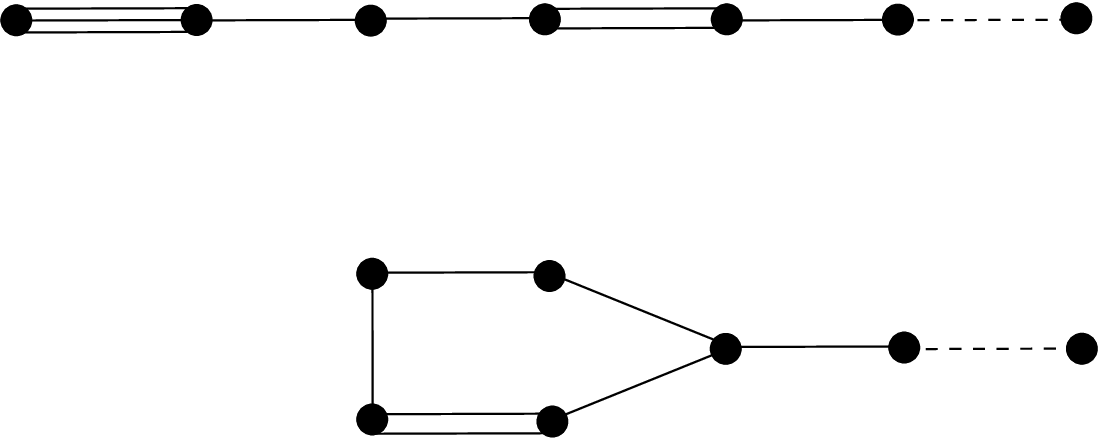,width=0.35\linewidth}\\
\hline
\end{tabular}
\caption{Straight hyperbolic Coxeter prisms}
\label{tab:all-prisms}
\end{table}

\section{Preliminaries}\label{sec:prelim}

\subsection{Hyperbolic lattices}

Let $\mathbb{R}^{d,1}$ be the real vector space $\mathbb{R}^{d+1}$ equipped with the standard quadratic form $f$ of signature $(d,1)$, namely, 
$$f(x)=-x_0^2+x_1^2+\dots+x_d^2.$$

The hyperboloid $\mathcal{H}=\{x \in \mathbb{R}^{d,1} \,|\, f(x) = -1 \}$ has two connected components 
$$
\mathcal{H}^+ = \{x \in \mathcal{H} \,|\, x_0 > 0\} \text{ and } \mathcal{H}^- = \{x \in \mathcal{H} \,|\, x_0 < 0\}.
$$

The $d$-dimensional hyperbolic space $\HH^d$ is the manifold $\mathcal{H}^+$ with the Riemannian metric $\rho$ induced by restricting $f$ to each tangent space $T_p(\mathcal{H}^+)$, $p \in \mathcal{H}^+$. This hyperbolic metric $\rho$ satisfies $\cosh \rho(x, y) = - (x, y)$, where $(x, y)$ is the scalar product in $\R^{d,1}$ associated to $f$.
The hyperbolic $d$-space $\HH^d$ is known to be the unique simply connected complete Riemannian $d$-manifold with constant sectional curvature $-1$.
{\em Hyperplanes} of $\HH^d$ are intersections of linear hyperplanes of $\R^{d,1}$ with $\mathcal{H}^+$, and are totally geodesic submanifolds of codimension $1$ in $\HH^d$. 

Let $\mathrm{O}_{d,1} = \mathbf{O}(f, \R)$ be the orthogonal group of the form $f$, and $\mathrm{O}'_{d,1} < \mathrm{O}_{d,1}$ be the subgroup (of index $2$) preserving $\mathcal{H}^+$. 
The group $\mathrm{O}'_{d,1}$ preserves the metric $\rho$ on $\HH^d$, and is in fact the full group $\mathrm{Isom}(\mathbb{H}^d)$ of isometries of the latter.

If $\Gamma < \mathrm{O}'_{d,1}$ is a lattice, i.e., if $\Gamma$ is a discrete subgroup of $\mathrm{O}'_{d,1}$ with a finite-volume fundamental domain in $\HH^d$, then the quotient $M=\mathbb{H}^d/\Gamma$ is a complete finite-volume {\itshape hyperbolic orbifold}. If $\Gamma$ is torsion-free, then $M$ is a complete finite-volume Riemannian manifold, and is called a {\itshape hyperbolic manifold}.

Denote $\mathrm{O}'_{d,1}$ by $G$, and suppose that $\G$ is an admissible (for $G$) algebraic $\mathbf{k}$-group, i.e. $\G(\R)^o$ is isomorphic to $G^o$ and $\G^\sigma(\R)$ is a compact group for any non-identity embedding $\sigma \colon \mathbf{k} \hookrightarrow \mathbb{R}$. Then any subgroup $\Gamma < G$ commensurable with the image in $G$ of $\G(\mathcal{O}_\mathbf{k})$ is an \textit{arithmetic  lattice} (in $G$) whose {\em adjoint trace field} is $\mathbf{k}$. 

Since $G$ also admits nonarithmetic lattices, we discuss some weaker notions of arithmeticity for lattices in $G$. Following Vinberg \cite{Vin67}, a lattice $\Gamma < G$ is called {\em quasi-arithmetic} if some finite-index subgroup of $\Gamma$ is contained in the image in $G$ of $\G(\mathbf{k})$, where $\G$ is some admissible algebraic $\mathbf{k}$-group, and is called  \textit{properly quasi-arithmetic} if $\Gamma$ is  quasi-arithmetic, but not arithmetic. 

We recall that the notion of quasi-arithmeticity is indeed broader than that of arithmeticity; as was mentioned in the introduction, the nonarithmetic closed hyperbolic manifolds constructed by Agol \cite{agol2006systoles}, Belolipetsky--Thomson~\cite{belolipetskythomson}, and Bergeron--Haglund--Wise \cite{bergeronhaglundwise} exist in all dimensions and, as observed by Thomson~\cite{Thomson}, are quasi-arithmetic. The first examples of properly quasi-arithmetic lattices in dimensions $3$, $4$, and $5$ were constructed by Vinberg \cite{Vin67} via reflection groups.

\subsection{Convex polyhedra and arithmetic properties of hyperbolic reflection groups}\label{sec:arith}
 
A {\em (hyperbolic) reflection group} is a discrete subgroup of $\mathrm{O}'_{d,1}$ generated by reflections in hyperplanes. A fundamental polyhedron (or a chamber) of any finite-covolume reflection group is a finite-volume {\em Coxeter polyhedron}, that is, a finite-volume (and thus finite-sided) convex polyhedron in which the dihedral angle between any two adjacent facets is an integral submultiple of $\pi$. Conversely, given a finite-volume Coxeter polyhedron $P \subset \HH^d$, the group generated by the reflections in all the supporting hyperplanes, or {\em walls}, of $P$ is a finite-covolume reflection group $\Gamma < \mathrm{O}'_{d,1}$ with fundamental chamber $P$. We thus frequently conflate finite-volume Coxeter polyhedra in $\HH^d$ with their corresponding lattices in $\mathrm{O}'_{d,1}$ (or their corresponding hyperbolic {\em reflective orbifolds}).

Let  $H_e = \{x \in \HH^d \mid (x,e)=0\}$ be a hyperplane in $\HH^d \subset \R^{d,1} $ whose linear span in $\R^{d,1}$ has normal vector $e \in \R^{d,1}$ with $(e,e)=1$, and $H_e^- = \{x \in \HH^d \mid (x,e) \le 0\}$ be the half-space associated with it. If
$P = \bigcap_{j=1}^N H_{e_j}^- $
is a finite-sided Coxeter polyhedron in $\HH^d$, 
then the matrix 
$$G(P) = \{g_{ij}\}^N_{i,j=1} = \{(e_i, e_j)\}^N_{i,j=1}$$ is its {\em Gram matrix}. We write $\mathbf{K}(P) = \Q\left(\{g_{ij}\}^N_{i,j=1}\right)$ and denote by $\mathbf{k}(P)$ the field generated by all possible cyclic products of the entries of $G(P)$; we call the field $\mathbf{k}(P)$ the {\em ground field} of $P$. For convenience, the set of all cyclic products of entries of a given matrix $A = (a_{ij})^N_{i,j=1}$, i.e., the set of all possible products of the form $a_{i_1 i_2} a_{i_2 i_3} \ldots a_{i_k i_1}$, will be denoted by $\Cyc(A)$. Thus, we have $\mathbf{k}(P) = \Q\left(\Cyc(G(P))\right) \subset \mathbf{K}(P)$.

The following statement follows from work of Vinberg, see~\cite[Section 4, Theorem 5]{Vin71} (see also a paper of Dotti \cite{dotti} with a detailed discussion around commensurability invariants of Coxeter groups and a variety of interesting examples). 

\begin{theorem}\label{Vin-field}
  Let $P$ be a Coxeter polyhedron in $\HH^d$ and suppose $\Gamma_P$ is the associated reflection group. If $\Gamma_P$ is Zariski-dense in $\Isom(\HH^d)$, then $\mathbf{k}(P)$ is the adjoint trace field of $\Gamma_P$, and thus is its commensurability invariant.
\end{theorem}

The following criterion allows us to determine if a given finite-covolume hyperbolic reflection group $\Gamma_P$ with fundamental polyhedron $P$ is arithmetic, quasi-arithmetic, or neither.

\begin{theorem}[Vinberg's arithmeticity criterion \cite{Vin67}]\label{V}
Let $\Gamma_P < \mathrm{O}
'_{d,1}$ be a reflection group with finite-volume fundamental polyhedron $P \subset \HH^d$. Then $\Gamma_P$ is arithmetic if and only if each of the following conditions holds:
\begin{itemize}
    \item[{\bf(V1)}] $\mathbf{K}(P)$ is a totally real algebraic number field;
    \item[{\bf(V2)}] for any embedding $\sigma \colon \mathbf{K}(P) \to \R$, such that 
    $\sigma\!\mid_{\mathbf{k}(P)} \ne \id$, the matrix $G^\sigma(P)$ is positive semi-definite;
    \item[{\bf(V3)}] $\Cyc(2 \cdot G(P)) \subset \OOO_{\mathbf{k}(P)}$.
\end{itemize}
The group $\Gamma_P$ is quasi-arithmetic if and only if it satisfies conditions {\em\textbf{(V1)}--\textbf{(V2)}}, but not necessarily {\em\textbf{(V3)}}.
\end{theorem}

\begin{remark}\label{rem:non-compact-QA}
    It is well-known that if a finite-volume Coxeter polyhedron $P$ is not compact, then quasi-arithmeticity of $\Gamma_P$ implies $\mathbf{k}(P) = \Q$. 
\end{remark}

\begin{remark}\label{rem:quasi-arithmetic}
Note that $2\cos \frac{\pi}{n}$ is always an algebraic integer. Thus, if there are no dashed edges in the Coxeter--Vinberg diagram of a finite-volume Coxeter polyhedron $P$, then condition  \textbf{(V3)} above automatically holds, and there is no distinction between arithmeticity and quasi-arithmeticity for the associated reflection group $\Gamma$. In particular, a triangle group acting on $\HH^2$ is quasi-arithmetic precisely when it is arithmetic.
\end{remark}

\begin{table}[]
\centering
{\renewcommand{\arraystretch}{1.77}
\begin{tabular}{|p{1.5cm}|p{4cm}|p{3cm}|}
\hline 
  $(k,l,m)$ & $a^2$ & Field $\mathbf{k}$ \\
\hline 
\hline 
  $(2,3,7)$ &  $\frac{3 \, \cos\left(\frac{1}{7} \, \pi\right)^{2} - 2}{4 \, \cos\left(\frac{1}{7} \, \pi\right)^{2} - 3}$  & $\mathbb{Q}(\cos(\frac{2\pi}{7}))$ \\
\hline
 $(2,3,8)$ & $ \frac{\sqrt{2} - 3}{2 \, {\left(\sqrt{2} - 2\right)}}$  & $\mathbb{Q}(\sqrt{2})$ \\
\hline 
  $(2,3,9)$ &  $\frac{3\cos(\frac{2\pi}{9})-1}{4\cos(\frac{2\pi}{9})-2}$ & $\mathbb{Q}(\cos(\frac{2\pi}{9}))$ \\ 
\hline 
  $(2,3,10)$ & $\frac{\sqrt{5}+7}{8} \,$ & $\mathbb{Q}(\sqrt{5})$ \\ 
\hline 
  $(2,3,14)$ & $\frac{3\cos(\frac{\pi}{7})-1}{4\cos(\frac{\pi}{7})-2}$ & $\mathbb{Q}(\cos(\frac{\pi}{7}))$  \\ 
\hline 
  $(2,4,5)$ & $\frac{\sqrt{5}+4}{4} \,$ & $\mathbb{Q}(\sqrt{5})$ \\ 
\hline 
  $(2,4,6)$ &$\frac{5}{4}$ & $\mathbb{Q}$\\ 
\hline 
  $(2,4,8)$ &  $\frac{\sqrt{2}+4}{4} \,$ & $\mathbb{Q}(\sqrt{2})$ \\
\hline 
  $(2,5,6)$ & $\frac{\sqrt{5}+7}{8} \,$ & $\mathbb{Q}(\sqrt{5})$\\ 
\hline 
  $(3,3,4)$ & $\frac{\sqrt{2}+1}{4}$  & $\mathbb{Q}(\sqrt{2})$ \\ 
\hline 
  $(3,3,5)$ & $\frac{\sqrt{5}+7}{8}$  & $\mathbb{Q}(\sqrt{5})$ \\ 
\hline 
  $(3,3,7)$  & $\frac{3 \, \cos\left(\frac{1}{7} \, \pi\right) - 1}{2 \, {\left(2 \, \cos\left(\frac{1}{7} \, \pi\right) - 1\right)}}$ & $\mathbb{Q}(\cos(\frac{\pi}{7}))$ \\   
\hline 
  $(3,3,9)$ & $\frac{3 \, \cos\left(\frac{1}{9} \, \pi\right) - 1}{2 \, {\left(2 \, \cos\left(\frac{1}{9} \, \pi\right) - 1\right)}}$ & $\mathbb{Q}(\cos(\frac{\pi}{9}))$ \\ 
\hline 
  $(4,4,3)$ &  $\frac{5}{4} $ & $\mathbb{Q}$ \\ 
\hline 
  $(4,4,4)$ &  $\frac{\sqrt{2}+3}{4} \,$ & $\mathbb{Q}(\sqrt{2})$ \\
\hline 
 $(4,5,2)$ & $\frac{\sqrt{5} + 1}{\sqrt{5} - 1}$ & $\mathbb{Q}(\sqrt{5})$ \\ 
\hline 
  $(4,5,4)$ & $\frac{3 \, {\left(\sqrt{5} + 2\right)}}{3 \, \sqrt{5} + 5}$ & $\mathbb{Q}(\sqrt{5})$ \\  
\hline 
  $(5,5,2)$ &  $\frac{\sqrt{5}}{\sqrt{5} - 1}$  & $\mathbb{Q}(\sqrt{5})$ \\ 
\hline 
  $(5,5,3)$ & $\frac{2 \, \sqrt{5} + 3}{2 \, {\left(\sqrt{5} + 1\right)}}$  & $\mathbb{Q}(\sqrt{5})$\\
\hline 
\end{tabular}}
\caption{Compact arithmetic prisms in $\HH^3$ of type $1$.}
\label{tab:type1-arith}
\end{table}

\section{Proof of Theorem~\ref{th:main}}\label{sec:main-proof}

Let us first discuss our general strategy. We will make use of Vinberg's arithmeticity criterion (see Theorem~\ref{V}) to pick quasi-arithmetic hyperbolic Coxeter prisms. We first notice that there are only finitely many such prisms in $\HH^4$ and $\HH^5$ while in $\HH^3$ we have several infinite families of them. Thus, in $\HH^4$ and $\HH^5$ it remains to apply Vinberg's arithmeticity criterion to a finite number of polyhedra.

\begin{table}
\centering
{\renewcommand{\arraystretch}{1.55}
\begin{tabular}{|p{1.5cm}|p{4cm}|p{3cm}| }
\hline
  $(k,l,m)$ & $a^2$ & Field $\mathbf{k}$\\ 
\hline 
\hline 
  $(2,3,12)$ & $\frac{\sqrt{3}+7}{8} \,$ & $\mathbb{Q}(\sqrt{3})$ \\ 
\hline 
  $(2,3,18)$ &  $\frac{3\cos(\frac{\pi}{9})-1}{4\cos(\frac{2\pi}{9})-2}$ & $\mathbb{Q}(\cos(\frac{\pi}{9}))$ \\ 
\hline 
  $(2,3,24)$ & $\frac{\sqrt{3}+\sqrt{2}+5}{8} \,$ & $\mathbb{Q}(\sqrt{3},\sqrt{2})$ \\ 
\hline 
  $(2,3,30)$ & $\frac{3\cos(\frac{\pi}{15})-1}{4\cos(\frac{\pi}{15})-2}$ & $\mathbb{Q}(\cos(\frac{\pi}{15}))$ \\ 
\hline 
  $(2,4,12)$ & $\frac{2\sqrt{3}+9}{12} \,$ & $\mathbb{Q}(\sqrt{3})$ \\ 
\hline 
  $(2,5,5)$ & $\frac{\sqrt{5}+8}{8} \,$ & $\mathbb{Q}(\sqrt{5})$ \\  
\hline 
  $(3,3,6)$ & $\frac{\sqrt{3}+7}{8}$  & $\mathbb{Q}(\sqrt{3})$\\
\hline
 $(3,3,12)$ & $\frac{\sqrt{2}+\sqrt{3}+5}{8} \,$ & $\mathbb{Q}(\sqrt{3},\sqrt{2})$  \\ 
\hline 
 $(3,3,12)$ & $\frac{\sqrt{2}+\sqrt{3}+5}{8} \,$ & $\mathbb{Q}(\sqrt{3},\sqrt{2})$  \\
\hline
$(3,4,3)$ & $\frac{3}{8} \, \sqrt{2} + 1$ &  $\mathbb{Q}(\sqrt{2})$ \\  
\hline 
  $(3,4,4)$ & $\frac{7}{6}$ & $\mathbb{Q}$ \\ 
\hline 
  $(3,4,6)$ & $\frac{\sqrt{6}+6}{8}$ & $\mathbb{Q}(\sqrt{2},\sqrt{3})$ \\ 
\hline 
  $(3,5,5)$ & $\frac{11 \, \sqrt{5} + 17}{12 \, {\left(\sqrt{5} + 1\right)}}$ & $\mathbb{Q}(\sqrt{5})$\\ 
\hline 
 $(4,4,6)$ &  $\frac{2\sqrt{3}+9}{12} \,$ & $\mathbb{Q}(\sqrt{3})$\\
 \hline 
 $(5,5,5)$ & $\frac{19 \, \sqrt{5} + 25}{20 \, {\left(\sqrt{5} + 1\right)}}$  & $\mathbb{Q}(\sqrt{5})$\\
\hline
\end{tabular}}
\caption{Compact properly quasi-arithmetic prisms in $\HH^3$ of type~$1$.}
\label{tab:type1-quasi}
\end{table}

To rule out quasi-arithmeticity of noncompact Coxeter prisms, we simply need to check if the ground field $\mathbf{k}(P)$ is $\Q$; see Remark~\ref{rem:non-compact-QA}. Obviously, this is not the case if the polyhedron $P$ has a dihedral angle of the form $\pi/m$ with $m > 6$ and $m=5$, as only $m = 2, 3, 4, 6$ gives rational $\cos^2 (\frac{\pi}{m})$ (see, for example, \cite[Corollary 3.12]{Niven}). In particular, this immediately implies that type-$7$ noncompact Coxeter prisms in $\HH^3$ are not quasi-arithmetic. Thus, in the noncompact case, even in $\HH^3$ it remains to check only a finite number of prisms.

Now we need to deal with several infinite families of compact Coxeter prisms in $\HH^3$. In this case we are going to apply a recent result of Bogachev and Kolpakov, see \cite{BK21}. It was shown in \cite{BK21} that a face of a quasi-arithmetic Coxeter polyhedron should also correspond to a quasi-arithmetic reflection group if it is itself a Coxeter polyhedron. For quasi-arithmetic straight Coxeter prisms this implies, in particular, that their triangle base (which is orthogonal to its neighbors) should also be quasi-arithmetic, and even arithmetic according to Remark~\ref{rem:quasi-arithmetic}. In fact, one could even apply even a more recent result in this situation: the above mentioned triangle base $F$ is a totally geodesic reflective suborbifold of an ambient prismatic Coxeter orbifold $P$, and by \cite[Theorem 1.7]{BBKS21}, if $P$ is quasi-arithmetic, then $F$ should be quasi-arithmetic as well.

That is, according to Remark~\ref{rem:quasi-arithmetic}, nonarithmeticity of the triangle reflection group associated to $F$ implies non-quasi-arithmeticity of $P$. By the result of Takeuchi \cite{Tak77}, only finitely many triangle groups in $\HH^2$ are arithmetic. This leaves us only finitely many compact Coxeter prisms in $\HH^3$ to be combed through and verify on quasi-arithmeticity following Vinberg's criterion (Theorem~\ref{V}). 

Let us now provide more details.

\subsection{Computing distances between disjoint facets}

By the theorem of Andreev \cite{An70}, all acute-angled prisms in $\HH^d$ have the following feature: there is exactly one pair of divergent facets (the bases of prisms). To use Vinberg's criterion we need to find this distance $\cosh d_P$ for every prism $P$, since it would give us the only one unknown entry of the Gram matrix of $P$. Any such Gram matrix $G(P)$ has signature $(d,1,1)$ which gives us the condition $\det G(P) = 0$. This is in fact a quadratic equation with respect to precisely one variable: $a = \cosh d_P$. For infinite families of compact Coxeter prisms (i.e. for types $1$--$3$) we collect the results of our computations in Table~\ref{tab:systoles} ($\cosh^2 d_m$ is presented for every fixed type and parameters $k,l$). In all other cases we have only a finite number of prisms for which this unknown entry need to be computed.

\subsection{Noncompact Coxeter prisms}

Recall that if $P$ is noncompact and quasi-arithmetic, then $\mathbf{k}(P) = \Q$ and therefore $P$ has no dihedral angles $\frac{\pi}{m}$ for $m = 5$ or $m > 6$. Let us also notice that $\cos \frac{\pi}{4} = \frac{1}{\sqrt{2}}$, $\cos \frac{\pi}{6} = \frac{\sqrt{3}}{2}$, and thus these angles may participate in cyclic products with another numbers only if the latter contain $\sqrt{2}$ or $\sqrt{3}$, respectively.

Noncompact prisms in $\HH^3$ have types $5$--$11$. For prisms of type $5$ it gives precisely two options, $m=4$ and $m=6$, corresponding to arithmetic and properly quasi-arithmetic reflection groups, respectively. The only prism of type $6$ has $m=6$: it is arithmetic. It was mentioned earlier that type-$7$ prisms are not quasi-arithmetic. For type-$8$ prisms, we consider the following non-trivial cyclic product $\cos \frac{\pi}{4} \cos \frac{\pi}{k} \cos \frac{\pi}{m}$:  if $k=3$ then $m=4$, since otherwise, this cyclic product belongs to some extension of $\Q$, and if $k=4$ then only $m=3$ or $m=6$ are possible. Thus, for the type $8$ only the following options are available to give a quasi-arithmetic reflection group: $(k,4,m) = (3,4,4), (4,4,3), (4,4,6)$. Similar considerations for the type $9$ leave only $m=6$ to be checked. The Gram matrix of type-$10$ prisms produce the cyclic product $\cos \frac{\pi}{k} \cos \frac{\pi}{l} \cos \frac{\pi}{3}$. If it belongs to $\Q$, then $k=l=4$ which gives an arithmetic group. The unique type-$11$ prism is also arithmetic. 

Noncompact prisms in $\HH^4$ have types $17$--$20$. The type-$17$ prism has $\sqrt{2}$ in the cyclic product, thus the Vinberg field is not $\Q$. The type-$18$ prism is quasi-arithmetic only for $k=4$. Finally, prisms of type $19$ and $20$ are not quasi-arithmetic, since they both have the angle $\pi/5$ contributing $\sqrt{5}$ to its Vinberg field.
Similar considerations prove non-quasi-arithmeticity of both noncompact prisms in $\HH^5$.

\subsection{Compact Coxeter prisms}

There are only finitely many compact Coxeter prisms in $\HH^4$ and $\HH^5$, and for them we directly apply Vinberg's arithmeticity criterion. They all turn out to be quasi-arithmetic. 

Now we need to rule quasi-arithmeticity of the three infinite families of compact Coxeter prisms in $\HH^3$: types $1$--$3$. As discussed above, any of these prisms has a totally geodesic triangle face, which would also be quasi-arithmetic if the ambient prism were. For these types $1$--$3$, the parameters $k,l$ of the associated triangle groups $(k,l,m)$ are bounded from above: $k,l \leqslant 5$ for type $1$ and $k,l \leqslant 3$ for types $2$--$3$, see Table~\ref{tab:all-prisms}.
Under these conditions and assuming $k \leqslant l$, the classification of Takeuchi \cite{Tak77} gives the following list of compact arithmetic triangles $(k,l,m)$ in $\HH^2$:

$(2,3,7), (2,3,8), (2,3,9), (2,3,10), (2,3,11), (2,3,12),(2,3,14), (2,3,16), (2,3,18)$, 

$(2,3,24), (2,3,30), (2,4,5),
(2,4,6), (2,4,7), (2,4,8), (2,4,10), (2,4,12), (2,4,18),$ 

$(2,5,5), (2,5,6), (2,5,8), (2,5,10), (2,5,20), (2,5,30)$, $(3, 3, 4), (3, 3, 5), (3, 3, 6),$

$(3, 3, 7), (3, 3, 8), (3, 3, 9), (3, 3, 12)$, $(3, 3, 15)$, 
$(3,4,4), (3,4,6)$, $(3,4,12), (3,5,5)$, 

$(4,4,4),(4,4,5), (4,4,6), (4,4,9)$, $(4,5,5)$,  $(5, 5, 5), (5, 5, 10), (5, 5, 15)$.

Compact arithmetic Coxeter prisms of type $1$ are collected in Table \ref{tab:type1-arith}, properly quasi-arithmetic ones of the same type, in Table \ref{tab:type1-quasi}. Compact arithmetic and properly quasi-arithmetic prisms of type $2$--$4$ are presented in Table~\ref{tab:type2-4-quasi}. In many non-quasi-arithmetic cases the field  $\mathbf{K}(P)$ generated by the Gram matrix entries is not totally real, while in other remaining cases the condition \textbf{(V2)} of Theorem \ref{V} does not hold. 
\qed

\begin{table}[]
\centering
{\renewcommand{\arraystretch}{1.5}
\begin{tabular}{|p{1.5cm}|p{3cm}|p{2cm}|p{1.6cm}|}
\hline 
$(k,l,m)$ & $a^2$  & Field $\mathbf{k}$  &A/PQA\\
\hline 
\hline
\multicolumn{4}{|c|}{Type 2} \\
\hline 
 $(2,3,7)$ & $\frac{2 \, \cos\left(\frac{\pi}{7}\, \right)^{2} - 1}{4 \, \cos\left(\frac{\pi}{7} \, \right)^{2} - 3}$ & $\mathbb{Q}(\cos(\frac{\pi}{7}))$&  PQA\\
\hline 
  $(2,3,8)$ & $\frac{\sqrt{2}+2}{2}$ & $\mathbb{Q}(\sqrt{2})$ &  A\\
\hline 
  $(2,3,10)$ & $\frac{\sqrt{5}+3}{4}$ & $\mathbb{Q}(\sqrt{5})$ & A\\
\hline 
  $(2,3,12)$ & $\frac{\sqrt{3}+3}{4}$ & $\mathbb{Q}(\sqrt{3})$ & A\\
\hline 
  $(2,3,18)$ &  $\frac{2 \, \cos\left(\frac{1}{18} \, \pi\right)^{2} - 1}{4 \, \cos\left(\frac{1}{18} \, \pi\right)^{2} - 3}$ & $\mathbb{Q}(\cos(\frac{\pi}{9}))$ & PQA\\
\hline 
  $(3,3,4)$ &$\frac{\sqrt{2}+2}{2}$  & $\mathbb{Q}(\sqrt{2})$ &  A\\
\hline 
 $(3,3,5)$ & $\frac{\sqrt{5}+3}{4}$  & $\mathbb{Q}(\sqrt{5})$ & A\\
\hline 
 $(3,3,6)$ & $\frac{\sqrt{3}+3}{4}$  & $\mathbb{Q}(\sqrt{3})$ &  A\\
\hline 
  $(3,3,9)$ & $\frac{\cos\left(\frac{\pi}{9} \, \right)}{2 \, \cos\left(\frac{\pi}{9} \, \right) - 1}$ & $\mathbb{Q}(\cos(\frac{\pi}{9}))$  & PQA\\
\hline
\multicolumn{4}{|c|}{Type 3} \\
\hline 
 $(3,3,5)$ & $\frac{1}{40} \, \sqrt{5} {\left(9 \sqrt{5} + 5\right)}$ & $\mathbb{Q}(\sqrt{5})$ &  PQA\\
\hline
\multicolumn{4}{|c|}{Type 4} \\
\hline
 $m=5$ & $\frac{\sqrt{5}+3}{2}$  & $\mathbb{Q}(\sqrt{5})$ & A \\ 
\hline
\end{tabular}}
\caption{Arithmetic (A) and properly quasi-arithmetic (PQA) compact Coxeter prisms in $\HH^3$ of types $2$, $3$, and $4$.}
\label{tab:type2-4-quasi}
\end{table}

\begin{table}[]
\centering
{\renewcommand{\arraystretch}{1.4}
\begin{tabular}{|p{4cm}|p{1cm}|p{2cm}|p{1.3cm}|}
\hline 
Type of $P$ and $(k,l,m)$ & $a^2$  & Field $\mathbf{k}$  & A/PQA\\
\hline
\hline
 Type $5$; $m=4$  & $\frac{3}{2}$ & $\mathbb{Q}$  & A\\
\hline 
 Type $5$; $m=6$ & $\frac{9}{8}$ & $\mathbb{Q}$ & PQA\\
\hline
Type $6$; $m = 6$  & $\frac{3}{2}$ & $\mathbb{Q}$  & A\\
\hline 
Type $8$; $m=3$  &  $\frac{3}{2}$ & $\mathbb{Q}$ & A\\
\hline 
Type $8$;  $m=4$   &  $\frac{4}{3}$ & $\mathbb{Q}$  & PQA\\
\hline
Type $10$; $k =l= 4$  & $2$ & $\mathbb{Q}$ & A\\
\hline 
Type $11$   & $\frac{3}{2}$ &  $\mathbb{Q}$ & A\\
\hline
\end{tabular}}
\caption{Arithmetic (A) and properly quasi-arithmetic (PQA) {\em noncompact} Coxeter prisms in $\HH^3$.}
\label{tab:types-5-11}
\end{table}

\begin{table}[]
\centering
{\renewcommand{\arraystretch}{1.47}
\begin{tabular}{|p{2.5cm}|p{3cm}|p{2cm}|p{1.5cm}|}
\hline 
Type of $P$ and $k$ & $a^2$  & Field $\mathbf{k}$  & A/PQA\\
\hline
\hline
\multicolumn{4}{|c|}{Types 12--16, compact prisms in $\HH^4$} \\ 
\hline
Type $12$ & $\frac{\sqrt{5}+3}{4}$  & $\mathbb{Q}(\sqrt{5})$  & A\\
\hline 
Type $13$  & $\frac{21+3\sqrt{5}}{22}$  & $\Q(\sqrt{5})$  & PQA\\
\hline 
Type $14$  & $\frac{\sqrt{5}+7}{8}$  & $\mathbb{Q}(\sqrt{5})$  & A\\
\hline 
Type $15$, $k=4$  & $\frac{6+2\sqrt{2}}{7}$  & $\mathbb{Q}(\sqrt{2})$   & PQA\\
\hline 
Type $15$, $k=5$ & $\frac{4\sqrt{5}+24}{31}$  & $\mathbb{Q}(\sqrt{5})$   & PQA\\
\hline 
Type $16$ & $\frac{\sqrt{5}+3}{2}$  & $\mathbb{Q}(\sqrt{5})$   & A\\
\hline
\multicolumn{4}{|c|}{Types 17--20, noncompact prisms in $\HH^4$} \\ 
\hline
Type $18$, $k=4$ & $\frac{8}{7}$ & $\mathbb{Q}$ &  PQA\\
\hline 
\multicolumn{4}{|c|}{Types 21--22, compact prisms in $\HH^5$} \\ 
\hline
Type $21$, $k=3$ & $\frac{\sqrt{5} + 7}{8}$ & $\mathbb{Q}(\sqrt{5})$  & A\\
\hline 
Type $21$, $k=4$ & $\frac{\sqrt{5} + 3}{4}$ & $\mathbb{Q}(\sqrt{5})$  & A\\
\hline 
Type $22$ & $\frac{\sqrt{2} + 3}{4}$ & $\mathbb{Q}(\sqrt{2})$  & A\\
\hline 
\multicolumn{4}{|c|}{Types 23--24, noncompact prisms in $\HH^5$: no one is QA} \\ 
 \hline
\end{tabular}}
\caption{Arithmetic (A) and properly quasi-arithmetic (PQA) Coxeter prisms in $\HH^4$ and $\HH^5$.}
\label{tab:other-types}
\end{table}

\begin{table}
\centering
{\renewcommand{\arraystretch}{1.5}
\begin{tabular}{|l|l|l|l|l|l|l|l|l|}
\hline
\multicolumn{3}{|c|}{Type 1} & \multicolumn{3}{|c|}{Type 2} & \multicolumn{3}{|c|}{Type 3}\\ 
\hline 
  $k$ & $l$ & $\cosh^2 d_m$ & $k$ & $l$ & $\cosh^2 d_m$  & $k$ & $l$ & $\cosh^2 d_m$ \\
\hline
\hline 
 $2$ & $3$  & $\frac{3\cos(\frac{2\pi}{m})-1}{4\cos(\frac{2\pi}{m})-2}$ & $2$ & $3$  & $\frac{2 \, \cos\left(\frac{\pi}{m}\right)^{2} - 1}{4 \, \cos\left(\frac{\pi}{m}\right)^{2} - 3}$ &   $2$ & $3$  & $\frac{-\left(\left(\sqrt{5}-5\right) \cos ^2\left(\frac{\pi }{m}\right)\right)+\sqrt{5}-3}{8 \cos ^2\left(\frac{\pi }{m}\right)-6}$   \\ 
 \hline 
 $2$ & $4$ & $\frac{3\cos(\frac{2\pi}{m})+1}{4\cos(\frac{2\pi}{m})}$  & $3$ &  $3$ & $\frac{\cos\left(\frac{\pi}{m}\right)}{2 \, \cos\left(\frac{\pi}{m}\right) - 1}$ & $3$ & $3$ & $\frac{\left(\sqrt{5}-5\right) \left(-\cos \left(\frac{\pi }{m}\right)\right)+\sqrt{5}-1}{8 \cos \left(\frac{\pi }{m}\right)-4}$ \\ 
\hline
 $2$ & $5$ & $\frac{3 \cos \left(\frac{2 \pi }{m}\right)+\sqrt{5}}{4 \cos \left(\frac{2 \pi }{m}\right)+\sqrt{5}-1}$  & &   &  & &   & \\ 
\hline
 $3$ & $3$ & $\frac{1-3 \cos \left(\frac{\pi }{m}\right)}{2-4 \cos \left(\frac{\pi }{m}\right)}$  & &   &  & &   &  \\ 
\hline
 $3$ & $4$ & $\frac{3e^2+2\sqrt{2}e}{4e^2+2\sqrt{2}e-1}$  & &   &  & &   &    \\ 
\hline
 $3$ & $5$ & $\frac{6e^2+2e+2\sqrt{5}e-1}{8e^2+2\sqrt{5}e+2e-3}$  & &   &  & &   &    \\ 
\hline
 $4$ & $4$ & $\frac{3 \, \cos\left(\frac{\pi}{m}\right) + 1}{4 \, \cos\left(\frac{\pi}{m}\right)}$  & &   &  & &   &  \\ 
\hline
 $4$ & $5$ & $\frac{\sqrt{5} + 3 \, \cos\left(\frac{\pi}{m}\right)}{\sqrt{5} + 4 \, \cos\left(\frac{\pi}{m}\right) - 1}$  & &   &  & &   &  \\ 
\hline
 $5$ & $5$ & $\frac{3 \cos \left(\frac{\pi }{m}\right)+\sqrt{5}}{4 \cos \left(\frac{\pi }{m}\right)+\sqrt{5}-1}$  & &   &  & &   &  \\ 
\hline
\end{tabular}}
\caption{Upper bounds on systoles of $\Gamma^{j}_{k,l,m}$. Here $e = \cos(\pi/m)$.}
\label{tab:systoles}
\end{table}

\section{Proof of Theorem~\ref{th:non-quasi}}\label{sec:non-quasi}

Let $P$ be a Coxeter prism glued from two straight prisms $P_1$ and $P_2$ along their common base $F$ that was orthogonal to its adjacent facets both in $P_1$ and $P_2$. Then it is clear that $F$ is also a totally geodesic reflective suborbifold of $P$. To rule out quasi-arithmeticity of $P$, one can again apply the above mentioned \cite[Theorem~1.7]{BBKS21}.

In our case, for fixed $k,l,m$, we obtain a Coxeter prism $P = P^{j,3}_{k,l,m}$, with $j = 1, 2$, by glueing two straight Coxeter prisms $P^{j}_{k,l,m}$ and $P^{3}_{k,l,m}$ along its common $(k,l,m)$ triangular face $F$. The Vinberg field $\mathbf{k}(F)$ is generated by $\cos (\frac{2\pi}{k})$, $\cos (\frac{2\pi}{l})$, $\cos (\frac{2\pi}{m})$, and $\cos (\frac{\pi}{k}) \cdot \cos (\frac{\pi}{l}) \cdot \cos (\frac{\pi}{m})$. If $m$ is coprime with $5$, then $\mathbf{k}(F)$ does not contain $\sqrt{5}$. On the other hand, the field $\mathbf{k}(P)$ definitely contains $\sqrt{5}$ since the dihedral angle $\pi/5$ appears in the Coxeter prism $P^{3}_{k,l,m}$. Thus, $\mathbf{k}(P) \not\subseteq \mathbf{k}(F)$, while it would be by \cite[Theorem 1.7]{BBKS21} were $P$ quasi-arithmetic.~\qed

\begin{remark}
    One may also observe that if $m$ is coprime with $5$ then $P^{j}_{k,l,m}$ and $P^{3}_{k,l,m}$ are not commensurable by the same argument via trace fields: $\mathbf{k}(P^{j}_{k,l,m}) \neq \mathbf{k}(P^{3}_{k,l,m})$ since the first field does not contain $\sqrt{5}$ while the second one does.
    Thus, the facet $F$ (the gluing locus) of $P = P^{j,3}_{k,l,m}$ is not an fc-subspace of the orbifold $P$ and therefore $P$ is not arithmetic, see \cite[Theorem 1.1]{BBKS21} (cf. also \cite[Theorem 1.1]{BGV23}).
\end{remark}

\section{Proof of Theorem~\ref{th:systoles}}\label{sec:systoles}

Given a hyperbolic lattice $\Gamma < \mathrm{Isom}(\mathbb{H}^d)$, the {\em systole} $\mathrm{sys}(\Gamma)$ of $\Gamma$ is the minimal translation length of a loxodromic element of $\Gamma$. The {\em systole} $\mathrm{sys}(M)$ of the complete finite-volume hyperbolic orbifold $M = \mathbb{H}^d / \Gamma$ is simply the systole of the lattice $\Gamma$. For complete finite-volume hyperbolic manifolds this definition gives precisely the well-known definition in the sense of shortest closed geodesics. For reflective orbifolds, it can be reformulated in the context of minimal closed billiard trajectories within the associated Coxeter polyhedra.

In Section~\ref{sec:main-proof}, we presented our metric computations for prisms which we used for Vinberg's arithmeticity criterion. In particular, we found the distances between opposite bases of prisms; see Table~\ref{tab:systoles}.

Let $P = \HH^d/\Gamma$ be a finite-volume hyperbolic reflective orbifold where $\Gamma$ is generated by reflections in the walls of the associated finite-volume Coxeter polyhedron $P$. It is clear that $\mathrm{sys}(P)$ is at most twice the minimal distance between disjoint facets of $P$. For a polyhedron $P^{j}_{k,l,m}$ denote such a distance by $d^{j}_{k,l,m}$.

It remains to observe that, given $j = 1, 2, 3$ and $k, l$, as in Table~\ref{tab:all-prisms} and Table~\ref{tab:systoles}, we have $\cosh^2 d^{j}_{k,l,m} \to 1$ as $m \to +\infty$, and therefore 
$$
\mathrm{sys}(P^{j}_{k,l,m}) \leq 2 d^{j}_{k,l,m} \xrightarrow[m \to \infty]{} 0. 
$$\qed

\begin{remark}
    Obviously, Theorem \ref{th:systoles} shows that there is no lower bound for systoles of compact reflective orbifolds in $\HH^3$. On the other hand, the recent paper \cite{B23} of the first author 
    demonstrates that apparently there may exist a universal upper bound for systoles of such orbifolds. More precisely, following the terminology of  \cite{B23}, if $P$ is a compact Coxeter polyhedron in $\HH^3$ with a small ridge such that the framing facets of the associated edge are divergent, then $\mathrm{sys}(P) < 2 \cdot \mathrm{arccosh}(5.75)$.
\end{remark}

\bibliography{biblio.bib}{}

\begin{thebibliography}{10}

\bibitem{agol2006systoles}
{\sc I.~Agol}, {\em Systoles of hyperbolic 4-manifolds}, arXiv preprint
  math/0612290,  (2006).

\bibitem{ABSW}
{\sc I.~Agol, M.~Belolipetsky, P.~Storm, and K.~Whyte}, {\em Finiteness of
  arithmetic hyperbolic reflection groups}, Groups Geom. Dyn., 2 (2008),
  pp.~481--498.

\bibitem{An70}
{\sc E.~M. Andreev}, {\em The intersection of the planes of the faces of
  polyhedra with sharp angles}, Mat. Zametki, 8 (1970), pp.~521--527.

\bibitem{BBKS21}
{\sc M.~Belolipetsky, N.~Bogachev, A.~Kolpakov, and L.~Slavich}, {\em Subspace
  stabilisers in hyperbolic lattices}, arXiv preprint arXiv:2105.06897,
  (2021).

\bibitem{belolipetskythomson}
{\sc M.~V. Belolipetsky and S.~A. Thomson}, {\em Systoles of hyperbolic
  manifolds}, Algebr. Geom. Topol., 11 (2011), pp.~1455--1469.

\bibitem{bergeronhaglundwise}
{\sc N.~Bergeron, F.~Haglund, and D.~T. Wise}, {\em Hyperplane sections in
  arithmetic hyperbolic manifolds}, J. Lond. Math. Soc. (2), 83 (2011),
  pp.~431--448.

\bibitem{B23}
{\sc N.~Bogachev}, {\em From geometry to arithmetic of compact hyperbolic
  {C}oxeter polytopes}, Transform. Groups, 28 (2023), pp.~77--105.

\bibitem{BD23}
{\sc N.~Bogachev and S.~Douba}, {\em Geometric and arithmetic properties of
  {L}\"obell orbifolds}, arXiv preprint arXiv:2304.12590,  (2023).

\bibitem{BDR23}
{\sc N.~Bogachev, S.~Douba, and J.~Raimbault}, {\em Infinitely many
  commensurability classes of compact {C}oxeter polyhedra in $\mathbb{H}^4$ and
  $\mathbb{H}^5$}, arXiv preprint arXiv:2309.07691,  (2023).

\bibitem{BGV23}
{\sc N.~Bogachev, D.~Guschin, and A.~Vesnin}, {\em Arithmeticity of ideal
  hyperbolic right-angled polyhedra and hyperbolic link complements}, arXiv
  preprint arXiv:2307.07000,  (2023).

\bibitem{BK21}
{\sc N.~Bogachev and A.~Kolpakov}, {\em On faces of quasi-arithmetic {C}oxeter
  polytopes}, Int. Math. Res. Not. IMRN,  (2021), pp.~3078--3096.

\bibitem{dotti}
{\sc E.~Dotti}, {\em On the commensurability of hyperbolic {C}oxeter groups},
  Manuscripta Math., 173 (2024), pp.~203--222.

\bibitem{DDK23}
{\sc E.~Dotti, S.~T. Drewitz, and R.~Kellerhals}, {\em Cusp density and
  commensurability of non-arithmetic hyperbolic {C}oxeter orbifolds}, Discrete
  Comput. Geom., 69 (2023), pp.~873--895.

\bibitem{dottikolpakov2022}
{\sc E.~Dotti and A.~Kolpakov}, {\em Infinitely many quasi-arithmetic maximal
  reflection groups}, Proc. Amer. Math. Soc., 150 (2022), pp.~4203--4211.

\bibitem{Douba}
{\sc S.~Douba}, {\em Systoles of hyperbolic hybrids}, arXiv preprint
  arXiv:2309.16051,  (2023).

\bibitem{Emery17}
{\sc V.~Emery}, {\em On volumes of quasi-arithmetic hyperbolic lattices},
  Selecta Math. (N.S.), 23 (2017), pp.~2849--2862.

\bibitem{Esselmann}
{\sc F.~Esselmann}, {\em \"{U}ber die maximale {D}imension von
  {L}orentz-{G}ittern mit coendlicher {S}piegelungsgruppe}, J. Number Theory,
  61 (1996), pp.~103--144.

\bibitem{GL14}
{\sc T.~Gelander and A.~Levit}, {\em Counting commensurability classes of
  hyperbolic manifolds}, Geom. Funct. Anal., 24 (2014), pp.~1431--1447.

\bibitem{GPS}
{\sc M.~Gromov and I.~I. Piatetski-Shapiro}, {\em Non-arithmetic groups in
  {Lobachevsky} spaces}, Publ. Math., Inst. Hautes {\'E}tud. Sci., 66 (1988),
  pp.~93--103.

\bibitem{Kap74}
{\sc I.~M. Kaplinskaja}, {\em The discrete groups that are generated by
  reflections in the faces of simplicial prisms in {L}oba\v{c}evski\u{\i}
  spaces}, Mat. Zametki, 15 (1974), pp.~159--164.

\bibitem{McR}
{\sc C.~Maclachlan and A.~W. Reid}, {\em The arithmetic of hyperbolic
  3-manifolds}, vol.~219 of Grad. Texts Math., New York, NY: Springer, 2003.

\bibitem{Mak66}
{\sc V.~S. Makarov}, {\em On a certain class of discrete groups of
  {L}oba\v{c}evski\u{\i}space having an infinite fundamental region of finite
  measure}, Dokl. Akad. Nauk SSSR, 167 (1966), pp.~30--33.

\bibitem{Mak68}
\leavevmode\vrule height 2pt depth -1.6pt width 23pt, {\em On {F}edorov's
  groups in four- and five-dimensional {L}obachevskij spaces}, Issled. po
  obshch. algebre. Kishinev, 1 (1968), pp.~120--129.

\bibitem{Mila18}
{\sc O.~Mila}, {\em Nonarithmetic hyperbolic manifolds and trace rings},
  Algebr. Geom. Topol., 18 (2018), pp.~4359--4373.

\bibitem{Nikulin}
{\sc V.~V. Nikulin}, {\em Finiteness of the number of arithmetic groups
  generated by reflections in {Lobachevsky} spaces}, Izv. Math., 71 (2007),
  pp.~53--56.

\bibitem{Niven}
{\sc I.~Niven}, {\em Irrational numbers}, vol.~No. 11 of The Carus Mathematical
  Monographs, Mathematical Association of America, ; distributed by John Wiley
  \& Sons, Inc., New York, 1956.

\bibitem{R13}
{\sc J.~Raimbault}, {\em A note on maximal lattice growth in {${\rm
  SO}(1,n)$}}, Int. Math. Res. Not. IMRN,  (2013), pp.~3722--3731.

\bibitem{Tak77}
{\sc K.~Takeuchi}, {\em Arithmetic triangle groups}, J. Math. Soc. Japan, 29
  (1977), pp.~91--106.

\bibitem{Thomson}
{\sc S.~Thomson}, {\em Quasi-arithmeticity of lattices in
  {{\(\mathrm{PO}(n,1)\)}}}, Geom. Dedicata, 180 (2016), pp.~85--94.

\bibitem{Vin67}
{\sc {\`E}.~B. Vinberg}, {\em Discrete groups generated by reflections in
  {L}oba\v{c}evski\u{\i}\ spaces}, Mat. Sb. (N.S.), 72 (114) (1967),
  pp.~471--488; correction, ibid. 73 (115) (1967), 303.

\bibitem{Vin71}
\leavevmode\vrule height 2pt depth -1.6pt width 23pt, {\em {Rings of definition
  of dense subgroups of semisimple linear groups}}, {Math. USSR, Izv.}, 5
  (1972), pp.~45--55.

\bibitem{vinberg84absence}
\leavevmode\vrule height 2pt depth -1.6pt width 23pt, {\em Absence of
  crystallographic groups of reflections in {L}obachevski\u{\i} spaces of large
  dimension}, Trudy Moskov. Mat. Obshch., 47 (1984), pp.~68--102, 246.

\end{thebibliography}
\bibliographystyle{siam}

\end{document}